\documentclass[12pt]{amsart}
\newtheorem{theorem}{Theorem}[section]
\newtheorem{lemma}[theorem]{Lemma}
\newtheorem{corollary}[theorem]{Corollary}
\newtheorem{proposition}[theorem]{Proposition}
\theoremstyle{definition}   
\newtheorem{definition}{Definition}
\newtheorem{example}[theorem]{Example}

\theoremstyle{remark}

\numberwithin{equation}{section}

\usepackage{times}
\usepackage{enumerate}
\usepackage{tfrupee}
\usepackage{tikz}
\usetikzlibrary{chains,fit}
\usetikzlibrary{shapes,snakes}
\usetikzlibrary{graphs}
\usepackage{graphicx,adjustbox}
\usepackage{amsmath,amssymb}
\usepackage[all,cmtip]{xy}
 \usepackage{calrsfs} 
\usepackage{tabularx}

\usepackage{amscd}
\usepackage{graphicx}
\usepackage[all]{xy}

\title[Numerical Semigroups with unique Ap\'{e}ry expansions]
{Numerical Semigroups with unique Ap\'{e}ry expansions}
\author{
Sudip Pandit
\and
Joydip Saha
\and
Indranath Sengupta
}
\date{}

\address{\small \rm  Discipline of Mathematics, IIT Gandhinagar, Palaj, Gandhinagar, 
Gujarat 382355, INDIA.}
\email{sudip.pandit@iitgn.ac.in}
\thanks{This work was done as an M.Sc. research project 
of the first author at IIT Gandhinagar, under 
the supervision of the third author, and in collaboration 
with the 
second author who was an SERB post-doctoral fellow at IIT Gandhinagar. 
The first author gratefully acknowledges the support of IIT Gandhinagar 
through the \textit{Summer 
research Internship Programme (SRIP)} in the year 2017, when 
this work began.}

\address{\small \rm  Stat-Math Unit, Indian Statistical Institute, 203 B.T. Road, Kolkata 700 108.} 
\email{saha.joydip56@gmail.com}
\thanks{The second author acknowledges the receipt of post-doctoral fellowship from 
SERB, Government of India.}

\address{\small \rm  Discipline of Mathematics, IIT Gandhinagar, Palaj, Gandhinagar, 
Gujarat 382355, INDIA.}
\email{indranathsg@iitgn.ac.in}
\thanks{The third author is the corresponding author; supported by the 
MATRICS research grant MTR/2018/000420, sponsored by the SERB, Government of India.}

\date{}



\allowdisplaybreaks

\begin{document}
\begin{abstract}
In this paper, we carry out a fairly comprehensive 
study of two special classes of numerical semigroups, 
one generated by the sequence of partial 
sums of an arithmetic progression and the other 
one generated by the partial sums of a geometric progression, 
in embedding dimension $4$. Both these 
classes have the common feature that they have unique expansions 
of the Ap\'{e}ry set elements. 
\end{abstract}

\maketitle

\section{Introduction}
The set $\Gamma$, subset of the set of nonnegative integers 
$\mathbb{N}$, is called a \textit{numerical semigroup} if it is 
closed under addition, contains zero and generates $\mathbb{Z}$ as 
a group. Every numerical semigroup $\Gamma$ satisfies the following two 
fundamental properties (see \cite{rs}): The complement $\mathbb{N}\setminus \Gamma$ is 
finite and $\Gamma$ has a unique minimal system of generators 
$a_{1} < \cdots < a_{n}$. The greatest integer not belonging to $\Gamma$, 
usually denoted by $F(\Gamma)$ is called the \textit{Frobenius number} 
of $\Gamma$. The integers $a_{1}$ and $n$, denoted by $m(\Gamma)$ and 
$e(\Gamma)$ respectively are known as the \textit{multiplicity} and the 
\textit{embedding dimension} of the semigroup $\Gamma$. The 
\textit{Ap\'{e}ry set} of $\Gamma$ with respect to a non-zero $a\in \Gamma$ is 
defined to be the set $\rm{Ap}(\Gamma,a)=\{s\in \Gamma\mid s-a\notin \Gamma\}$. 
In this paper, we study a class of numerical semigroups, which are special 
in the sense that they have the uniqueness of representation of each element 
in the Ap\'{e}ry set $\rm{Ap}(\Gamma,a)$. Given positive 
integers $a_{1} < \cdots < a_{n}$, every numerical semigroup ring 
$k[\Gamma] = k[t^{a_{1}}, \ldots , t^{a_{n}}]$ is the coordinate 
ring of an affine monomial curve given by the monomial parametrization 
$\nu : k[x_{1}, \ldots, x_{n}]\longrightarrow k[t]$, such that 
$\nu(x_{i}) = t^{a_{i}}$, $1\leq i\leq p$. The ideal $\ker(\nu)=\mathfrak{p}$ 
is the defining ideal of the parametrized monomial curve, which is graded 
with respect to the weighted gradation. 
\medskip

It is known that uniqueness of 
representations of the Ap\'{e}ry set elements of a numerical semigroup is 
actually quite helpful; see \cite{rs1}, \cite{mss}. Let us call these as 
numerical semigroups with unique Ap\'{e}ry expansions. 
One requires the Ap\'{e}ry table in order to 
understand the tangent cone, which is quite hard to compute in general. However, 
uniqueness of expressions of the Ap\'{e}ry set elements makes it easier. 
In this paper, we will be presenting two classes of numerical semigroups with this property. 
In fact, we have stumbled upon these classes while looking for a large 
class of numerical semigroups with this property, especially from the standpoint 
of computing tangent cones. 
\medskip

Let $f(x),g(x)\in\mathbb{Q}[x]$ such that  
$f(\mathbb{N})\subset \mathbb{N}$, $g(\mathbb{N})\subset \mathbb{N}$ and 
both are increasing, so called increasing numerical polynomials. In this paper 
we study numerical semigroups minimally generated by integers of the form 
$\{a, g(i)a+f(i)d\mid \gcd(a,d)=1, 1\leq i\leq n\}$. Some 
of the interesting classes of numerical semigroups that have been studied fall 
under this general class. For example, $g(x)=1$, $f(x)=x$ gives an 
arithmetic sequence (see \cite{patil}, \cite{pss}) and $g(x)$ = constant, 
$f(x)=x$ gives a generalized arithmetic sequence (see \cite{matthews}). 
We study the following two cases:
\medskip

\begin{enumerate}
\item When $g(x)=x+1$, $f(x)=\dfrac{x(x+1)}{2}$ and $n=3$; we denote the semigroup 
by $\Gamma_{4}$. A complete study has been carried out in sections 2 through 6 in the 
following oder - Ap\'{e}ry set, the pseudo Frobenius numbers, the 
defining ideal, syzygies and finally the Ap\'{e}ry table and the tangent cone. All 
these are known to be extremely hard to compute in general. We have used the computer 
algebra system \cite{GAP4} to form initial guesses for many of the theorems that we 
have proved here.  
\medskip

\item We take $g(x)$ to be a constant numerical function, $f(x) = r^{x}$ and 
define the numerical semigroup $\mathfrak{S}_{n+2}$, in section 7. We 
compute the the Ap\'{e}ry table and the tangent cone for $\mathfrak{S}_{n+2}$.
\end{enumerate}
\bigskip

\section{Ap\'{e}ry set of $\Gamma_{4}$}
We now consider the numerical semigroup generated 
by the positive integers $s_{1},\ldots, s_{4}$, where $d>0$ and $a>0$ are integers with 
$\gcd(a,d)=1$, and $s_{n}=\frac{n}{2}[2a+(n-1)d]$, for $1\leq n\leq 4$. We denote this numerical semigroup 
by $\Gamma_{4}$, the semigroup ring by $k[\Gamma_{4}]$ and the defining ideal 
by $\mathfrak{p}_{4}$. We will see in the next Theorem that we need to impose some 
bounds on $a$ so that $\{s_{1},\ldots, s_{4}\}$ is a minimal generating 
set for the numerical semigroup $\Gamma_{4}$. 
\medskip

\begin{theorem}\label{min4}
Let $d>0$ and $a\geq 7$ be integers with $\gcd(a,d)=1$. Let 
$s_{n}=\frac{n}{2}[2a+(n-1)d]$, for $1\leq n\leq 4$. The set 
$T_{4}=\{s_{1},\ldots, s_{4}\}$ is a minimal generating 
set for the numerical semigroup $\Gamma_{4} = \langle s_{1},\ldots, s_{4} \rangle$.
\end{theorem}

\proof 
Suppose $s_{3}=m_{1}s_{1}+m_{2}s_{2}$, for some $m_{1},m_{2}\geq 0$. We get,
\begin{eqnarray}\label{eq1}
(m_{1}+2m_{2}-3)a &=& (3-m_{2})d.
\end{eqnarray}
Since $\gcd(a,d)=1$, we get $a\mid (3-m_{2})$. 
\medskip

If $m_{2}\leq 3$, then $3=m_{2}+ka$, for some $k\geq 0$, and we get 
$m_{2}=3$, since $a\geq 7$. Therefore, $m_{1}+3 =0$ (using equation \ref{eq1}) - 
a contradiction. If $m_{2}> 3$, then the L.H.S. of equation 
\ref{eq1} is positive whereas the R.H.S. is negative - a contradiction.
\medskip

Suppose $s_{4}=m_{1}s_{1}+m_{2}s_{2}+m_{3}s_{3}$, for some $m_{1},m_{2},m_{3}\geq 0$.
We get,
\begin{eqnarray}\label{eq2}
(m_{1}+2m_{2}+3m_{3}-4)a &=& (6-m_{2}-3m_{3})d.
\end{eqnarray}
Therefore $a\mid (6-m_{2}-3m_{3})$, since $\gcd(a,d)=1$.  
\medskip

If $m_{2}+3m_{3}\leq 6$, then $6=m_{2}+3m_{3}+ka$, for some $k\geq 0$. Therefore, 
we get $6=m_{2}+3m_{3}$, since $a\geq 7$. Possible solutions for $(m_{2}, m_{3})$ 
are $(0,2)$, $(3,1)$, $(6,0)$. Substituting these values of 
$m_{2},m_{3}$ in the equation \ref{eq2}, we get $m_{1}<0$ - a contradiction.
\medskip

If $m_{2}+3m_{3}> 6$, then the L.H.S. of the equation \ref{eq2} is positive, 
whereas the R.H.S. is negative - a contradiction. \qed
\medskip

\begin{theorem}\label{apery}
Let $a\geq 7$. For each $1\leq i\leq a-1$, let $i=6\mu_{i}+q_{i}$, 
such that $0\leq q_{i}< 6$. For each $ 1\leq i\leq a-1 $, we define $\nu_{i},\xi_{i}$ 
as follows;
\begin{enumerate}
\item[(i)] $(\nu_{i},\xi_{i})=(1,q_{i}-3)$, if $q_{i}\geq 3$;
\item[(ii)] $(\nu_{i},\xi_{i})=(0,q_{i})$, if $q_{i}< 3$.
\end{enumerate}
Let $\mathrm{Ap}(\Gamma_{4},a)$ denote the Ap\'{e}ry set of $\Gamma_{4}$, 
with respect to the element $a$. Then 
$\mathrm{Ap}(\Gamma_{4},a)=\{(4\mu_{i}+3\nu_{i}+2\xi_{i})a+id \mid 1\leq i\leq a-1 \}\cup \{0\}$.
\end{theorem}

\proof Let $T=\{(4\mu_{i}+3\nu_{i}+2\xi_{i})a+id\mid 1\leq i\leq a-1 \}$. We notice that 
$i=6\mu_{i}+3\nu_{i}+\xi_{i}$, therefore for $1\leq i\leq a-1$, we have
$$(4\mu_{i}+3\nu_{i}+2\xi_{i})a+id=\mu_{i}(4a+6d)+\nu_{i}(3a+3d)+\xi_{i}(2a+d)\in \Gamma_{4}.$$ 
Hence $T\subset \Gamma_{4}$. Let $s\in \mathrm{Ap}(\Gamma_{4},a)\setminus\{0\}$, with 
$s\equiv id(\mbox{mod} \, a)$. 
Suppose
\begin{align*}
s&=c_{1}(2a+d)+c_{2}(3a+3d)+c_{3}(4a+6d)\\
&=(2c_{1}+3c_{2}+4c_{3})a+(c_{1}+3c_{2}+6c_{3})d, 
\end{align*}
then $(c_{1}+3c_{2}+6c_{3})\equiv i(\mbox{mod} \, a)$, as $\gcd(a,d)=1$. Therefore 
\begin{equation}\label{eqnap}
c_{1}+3c_{2}+6c_{3}= i +ka=6\mu_{i}+q_{i}+ka,
\end{equation}
for some $k\geq 0$. It is enough to show that, $4c_{3}+3c_{2}+2c_{1}\geq 4\mu_{i}+3\nu_{i}+2\xi_{i}$.
Suppose $4c_{3}+3c_{2}+2c_{1}< 4\mu_{i}+3\nu_{i}+2\xi_{i}$, then from \ref{eqnap}, 
substituting $\mu_{i}$, we have 
 \begin{equation}\label{ineq}
 6\xi_{i}+9\nu_{i}-2q_{i}>4c_{1}+3c_{2}+2ka .
 \end{equation}
We consider the following cases:
\medskip
 
\textbf{Case A.} If $q_{i}=0$, then $(\nu_{i},\xi_{i})=(0,0)$, and from \ref{ineq} we get 
$0>4c_{1}+3c_{2}+2ka$, which is impossible.
\medskip
 
\textbf{Case B.} If $q_{i}=1$, then $(\nu_{i},\xi_{i})=(0,1)$ and 
from \ref{ineq} we get $4>4c_{1}+3c_{2}+2ka$. Therefore $k=0$ and 
$(c_{1},c_{2})\in\{(0,0),(0,1)\}$. Putting values of $c_{1}$ and 
$c_{2}$ in equation \ref{eqnap}, we get the following:
$$
\begin{cases}
6c_{3}=6\mu_{i}+1 & \mbox{if} \quad (c_{1},c_{2})=(0,0),\\
6c_{3}=6\mu_{i}-2 & \mbox{if} \quad (c_{1},c_{2})=(0,1).
\end{cases}
$$
All lead to contradictions.
\medskip
 
\textbf{Case C.} If $q_{i}=2$, then $(\nu_{i},\xi_{i})=(0,2)$, and 
from \ref{ineq} we get $8>4c_{1}+3c_{2}+2ka$. 
Therefore $k=0$ and $(c_{1},c_{2})\in\{(0,0),(1,0),(1,1),(0,1),(0,2)\}$. 
Putting values of $c_{1}$ and $c_{2}$ in equation \ref{eqnap}, we get 
the following:
$$
\begin{cases}
6c_{3}=6\mu_{i}+2 & \mbox{if} \quad (c_{1},c_{2})=(0,0),\\
6c_{3}=6\mu_{i}+1 & \mbox{if} \quad (c_{1},c_{2})=(1,0),\\
6c_{3}=6\mu_{i}-2 & \mbox{if} \quad (c_{1},c_{2})=(1,1),\\
6c_{3}=6\mu_{i}-1 & \mbox{if} \quad (c_{1},c_{2})=(0,1),\\
6c_{3}=6\mu_{i}-4 & \mbox{if} \quad (c_{1},c_{2})=(0,2).
\end{cases}
$$
All lead to contradictions.
\medskip
 
\textbf{Case D.} If $q_{i}=3$, then $(\nu_{i},\xi_{i})=(1,0)$, and 
from \ref{ineq} we get $3>4c_{1}+3c_{2}+2ka$. Therefore $k=0$ and 
$(c_{1},c_{2})=(0,0)$. Then from equation \ref{eqnap}, we get $6c_{3}=6\mu_{i}+3$; 
which is not possible.
\medskip
 
\textbf{Case E.} If $q_{i}=4$, then $(\nu_{i},\xi_{i})=(1,1)$, and 
from \ref{ineq} we get $7>4c_{1}+3c_{2}+2ka$. Therefore $k=0$ and 
$(c_{1},c_{2})\in\{(0,0),(1,0),(0,1),(0,2)\}$. Putting values of 
$c_{1}$ and $c_{2}$ in equation \ref{eqnap}, we get the following:
$$
\begin{cases}
6c_{3}=6\mu_{i}+4 & \mbox{if} \quad (c_{1},c_{2})=(0,0),\\
6c_{3}=6\mu_{i}+3 & \mbox{if} \quad (c_{1},c_{2})=(1,0),\\
6c_{3}=6\mu_{i}+1 & \mbox{if} \quad (c_{1},c_{2})=(0,1),\\
6c_{3}=6\mu_{i}-2 & \mbox{if} \quad (c_{1},c_{2})=(0,2).
\end{cases}
$$
All lead to contradictions.
\medskip
 
\textbf{Case F.} If $q_{i}=5$, then $(\nu_{i},\xi_{i})=(1,2)$, from \ref{ineq} we get $11>4c_{1}+3c_{2}+2ka$. 
Therefore $k=0$ and 
$$(c_{1},c_{2})\in\{(0,0),(0,1),(0,2),(0,3),(1,0),(1,1),(1,2),(2,0)\}.$$ 
Putting 
values of $c_{1}$ and $c_{2}$ in equation \ref{eqnap}, we get the following:
$$
\begin{cases}
6c_{3}=6\mu_{i}+5 & \mbox{if} \quad (c_{1},c_{2})=(0,0),\\
6c_{3}=6\mu_{i}+2 & \mbox{if} \quad (c_{1},c_{2})=(0,1),\\
6c_{3}=6\mu_{i}-1 & \mbox{if} \quad (c_{1},c_{2})=(0,2),\\
6c_{3}=6\mu_{i}-4 & \mbox{if} \quad (c_{1},c_{2})=(0,3),\\
6c_{3}=6\mu_{i}+4 & \mbox{if} \quad (c_{1},c_{2})=(1,0),\\
6c_{3}=6\mu_{i}+1 & \mbox{if} \quad (c_{1},c_{2})=(1,1),\\
6c_{3}=6\mu_{i}-2& \mbox{if} \quad (c_{1},c_{2})=(1,2),\\
6c_{3}=6\mu_{i}+3 & \mbox{if} \quad (c_{1},c_{2})=(2,0).
\end{cases}
$$
All lead to contradictions. \qed
\medskip
 
\noindent An example has been discussed in \ref{atexample}.
\bigskip

\section{Pseudo Frobenius numbers and type of $\Gamma_{4}$}
\begin{definition}
Let $\Gamma$ be a numerical semigroup, we say thet $x\in\mathbb{Z}$ is a 
\textit{pseudo-Frobenius number} if $x\notin \Gamma$ and $x+s\in \Gamma$ 
for all $s\in \Gamma\setminus \{0\}$. We denote by $\mathbf{PF}(\Gamma)$ 
the set of  pseudo-Frobenius numbers of $\Gamma$. The cardinality of 
$\mathbf{PF}(\Gamma)$ is denoted by $t(\Gamma)$ and we call it the 
\textit{type} of $\Gamma$.
\end{definition}

Let $a,b\in \mathbb{Z}$. We define $\leq_{\Gamma}$ as $a\leq_{\Gamma} b$ 
if $b-a\in \Gamma$. This order relation defines a poset structure on 
$\mathbb{Z}$.

\begin{theorem}\label{max}
Let $\Gamma$ be a numerical semigroup and $a\in\Gamma\setminus \{0\}$. Then 
$\mathbf{PF}(\Gamma)=\{w-a\mid w\in \,\mathrm{Maximals}_{\leq_{\Gamma}}Ap(\Gamma, a)\}$. 
\end{theorem}
\proof See Proposition 8 in \cite{as}.\qed
\medskip

Let $\omega(i)=(4\mu_{i}+3\nu_{i}+2\xi_{i})a+id$, for $1\leq i\leq a-1$. 
Therefore, $\mathrm{Ap}(\Gamma_{4},a)=\{\omega(i)\mid 1\leq i\leq a-1 \}$.
\medskip

\begin{theorem}\label{min5}
Let  $a\geq 7$ and $d$ be two integers such that $\gcd(a,d)=1$. Suppose 
$\Gamma_{4}=\langle s_{1},\ldots, s_{4} \rangle$, where $s_{n}=\frac{n}{2}[2a+(n-1)d]$, 
for $1\leq n\leq 4$. Let $\mathbf{PF}(\Gamma_{4})$ be the set of pseudo Frobenius numbers 
of the numerical semigroup $\Gamma_{4}$. Write $a=6m+q$, $0\leq q\leq 5$; then 
\begin{eqnarray*}
\mathbf{PF}(\Gamma_{4}) & = & 
\begin{cases}
\{\omega(a-1)\}, \quad \mbox{if} \quad q=0;\\
\{\omega(a-1),\omega(a-2)\}, \quad \mbox{if} \quad q=1;\\
\{\omega(a-1),\omega(a-3)\}, \quad \mbox{if} \quad q=2;\\
\{\omega(a-1),\omega(a-4)\}, \quad \mbox{if} \quad q=3;\\
\{\omega(a-1),\omega(a-2),\omega(a-5)\}, \quad \mbox{if} \quad q=4;\\
\{\omega(a-1),\omega(a-3),\omega(a-6)\}, \quad \mbox{if} \quad q=5.
\end{cases}
\end{eqnarray*}
\end{theorem}

\proof  We first note that $\omega(i+6)- \omega(i)=4a+6d\in\Gamma_{4}$, 
for $0\leq i<a-6$. Therefore, $\omega(i)\leq_{\Gamma_{4}} \omega(i+6)$, 
$0\leq i<a-6$. Hence, by Theorem \ref{max}, $\mathbf{PF}(\Gamma_{4})\subset \{\omega(a-i)\mid 1\leq i\leq 6\}$. Proof of the theorem follows easily by checking each case.\qed
\medskip

\begin{corollary} Let $\mathrm{Der}_{k}(\Gamma_{4})$ be the set of $k$-derivations of $k[[t^{a},t^{2a+d},t^{3a+3d},t^{4a+6d}]]$, then
$$ \mathrm{Der}_{k}(\Gamma_{4})=\{t^{\alpha+1}\mid \alpha\in \mathbf{PF}(\Gamma_{4}) \}.$$
\end{corollary}

\proof Follows from \cite{kraft}, page 875.\qed 
\medskip

\begin{corollary}
Let $F(\Gamma_{4})$ be the Frobenius number of $\Gamma_{4}$. Then
\begin{enumerate}[(i)]
\item $F(\Gamma_{4})=\omega(a-1)$, \quad if $q\in\{0,3,5\}$;
\medskip

\item If $q=1$, then\\
$F(\Gamma_{4})=
\begin{cases}
\omega(a-2) & \mbox{if} \quad 3a>d;\\
\omega(a-1) & \mbox{otherwise}.
\end{cases}$
\medskip

\item If $q=2$ then\\
$F(\Gamma_{4})=
\begin{cases}
\omega(a-3) & \mbox{if} \quad a>2;\\
\omega(a-1) & \mbox{otherwise}.
\end{cases}$
\medskip

\item If $q=4$ then\\
$F(\Gamma_{4})=
\begin{cases}
\omega(a-2) & \mbox{if} \quad a>d;\\
\omega(a-1) & \mbox{otherwise}.
\end{cases}$
\end{enumerate}
\end{corollary}

\proof One can easily find the maximum element from \ref{min5}. \qed
\bigskip

\section{Minimal generating set for the defining ideal $\mathfrak{p}_{4}$}
Let us begin with the following theorem from \cite{g}, which helps us compute 
a minimal generating set for the defining ideal of a monomial curve.
\medskip

\begin{theorem}\label{gastinger}
 Let $A = k[x_{1},\ldots,x_{n}]$ be a polynomial ring, $I\subset A$ the defining
ideal of a monomial curve defined by natural numbers $a_{1},\ldots,a_{n}$, 
whose greatest common divisor is $1$.  Let $J \subset I$ be a subideal. 
Then $J = I$ if and only if $\mathrm{dim}_{k} A/\langle J + (x_{i}) \rangle =a_{i}$
for some $i$. (Note that the above conditions are also equivalent to 
$\mathrm{dim}_{k} A/\langle J + (x_{i}) \rangle =a_{i}$ for any $i$.)
\end{theorem}  

\proof See \cite{g}.\qed
\medskip
 
\begin{lemma}\label{equal}
Let $A = k[x_{1},\ldots,x_{n}]$ be a polynomial ring. For a monomial ideal $J$ of $A$, we write 
the unique minimal generating set of $J$ as $G(J)$. Let $I=\langle f_{1},\ldots f_{k}\rangle$ and 
$I_{i}=\langle f_{1},\ldots,\hat{f}_{i},\ldots f_{k}\rangle$, $1\leq i\leq k$. 
Suppose that, with respect to some monomial order on $A$, $\{ \mathrm{LT}(f_{1}),\ldots \mathrm{LT}(f_{k})\} \subset G(\mathrm{LT}(I)) $ and 
$G(\mathrm{LT}(I_{i})) \subset G(\mathrm{LT}(I))\setminus \{\mathrm{LT}(f_{i})\}$ for all $1\leq i\leq k$. 
Then $I$ is minimally generated by $\{f_{1},\ldots f_{k}\}$.
\end{lemma}

\proof Suppose $I$ is not minimally generated by $\{f_{1},\ldots f_{k}\}$. Then there is a 
polynomial $f_{i}$ such that $f_{i}\in I_{i}$. Therefore there is a monomial 
$m\in G(\mathrm{LT}(I_{i}))$, such that $m\mid \mathrm{LT}(f_{i})$. But $m$ and 
$\mathrm{LT}(f_{i})$ are distinct elements of $G(\mathrm{LT}(I))$, which gives a 
contradiction.  \qed
\medskip
 
\noindent\textbf{Notations.} We now introduce some notations specific to the polynomial ring with 
$4$ variables $A=k[x_{1},x_{2},x_{3},x_{4}]$. Let $m$ and $d$ be fixed positive integers. 
We define the subsets 
$H_{0}, H_{1}, H_{2}, H_{3}, H_{4}, H_{5}$ of $A$ as follows:
\medskip

\begin{enumerate}
\item[(i)] $H_{0}=\{x_{3}^{2}-x_{1}^{2}x_{4},x_{2}^{3}-x_{1}^{3}x_{3},x_{1}^{4m+d}-x_{4}^{m}\}$.
\medskip

\item[(ii)] 
\begin{itemize}
\item[(a)]
$H_{1}=\{x_{3}^{2}-x_{1}^{2}x_{4},x_{2}^{3}-x_{1}^{3}x_{3},x_{1}^{7}-x_{2}x_{4},x_{1}^{4}x_{2}^{2}-x_{3}x_{4},x_{1}^{2}x_{2}^{2}x_{3}-x_{4}^{2} \} $, if $m=d=1$.
\medskip

\item[(b)] 
$H_{1}=\{x_{3}^{2}-x_{1}^{2}x_{4},x_{2}^{3}-x_{1}^{3}x_{3},x_{1}^{(4m+d-6)}x_{2}^{5}-x_{4}^{m+1},x_{1}^{(4m+d-1)}x_{2}^{2}-x_{3}x_{4}^{m},x_{1}^{(4m+d+2)}-x_{2}x_{4}^{m}\}$, otherwise.
\end{itemize}
\medskip

\item[(iii)] $H_{2}=\{x_{3}^{2}-x_{1}^{2}x_{4},x_{2}^{3}-x_{1}^{3}x_{3},x_{1}^{(4m+d-4)}x_{2}^{4}-x_{4}^{m+1},x_{1}^{(4m+d+1)}x_{2}-x_{3}x_{4}^{m},x_{1}^{(4m+d+4)}-x_{2}^{2}x_{4}^{m}\}$.
\medskip

\item[(iv)] $H_{3}=\{x_{3}^{2}-x_{1}^{2}x_{4},x_{2}^{3}-x_{1}^{3}x_{3},x_{1}^{(4m+d-2)}x_{2}^{3}-x_{4}^{m+1},x_{1}^{(4m+d+3)}-x_{3}x_{4}^{m}\}$.
\medskip

\item[(v)] $H_{4}=\{x_{3}^{2}-x_{1}^{2}x_{4},x_{2}^{3}-x_{1}^{3}x_{3},x_{1}^{(4m+d)}x_{2}^{2}-x_{4}^{m+1},x_{1}^{(4m+d+5)}-x_{2}x_{3}x_{4}^{m}\}$.
\medskip

\item[(vi)] $H_{5}=\{x_{3}^{2}-x_{1}^{2}x_{4},x_{2}^{3}-x_{1}^{3}x_{3},x_{1}^{(4m+d+2)}x_{2}-x_{4}^{m+1},x_{1}^{(4m+d+7)}-x_{2}^{2}x_{3}x_{4}^{m}\}$.
\end{enumerate}
\medskip

\begin{theorem}\label{mingen}
Suppose $a=6m+q$, where $0\leq q\leq5$. Then $H_{q}$ is a minimal generating set for the 
ideal $\mathcal{P}_{4}$.
\end{theorem}

\proof We now use Theorem \ref{gastinger} to show that $\mathrm{dim}_{k}(A/\langle H_{q},x_{1}\rangle)=a$. 
Let $B=k[x_{2},x_{3},x_{4}]$ and $H^{'}_{q}=\langle H_{q},x_{1}\rangle/\langle x_{1}\rangle\subset B$. Therefore we need to show that  $\mathrm{dim}_{k}(B/\langle H_{q}^{'}\rangle)=a$. Let $\kappa_{q}$ be the dimension of the vector space $B/\langle H_{q}^{'}\rangle$.  We define 
$$\mathfrak{B}=\{x_{2}^{i}x_{3}^{j}x_{4}^{k}\mid 0\leq i\leq 2,0\leq j\leq 1,0\leq k\leq m \}.$$ 
and show that image of the set  $\mathfrak{B}\setminus\mathfrak{B}_{q}$ forms a basis of the 
vector space $B/\langle H_{q}^{'}\rangle$, through the following cases:

\begin{enumerate}[(A)]
\item We have $H_{0}^{'}=\{x_{3}^{2},x_{2}^{3},x_{4}^{m}\}$ and $\mathfrak{B}_{0}=\{x_{4}^{m},x_{2}x_{4}^{m}, x_{2}^{2}x_{4}^{m},x_{3}x_{4}^{m},x_{2}x_{3}x_{4}^{m}, x_{2}^{2}x_{3}x_{4}^{m}\} $. Hence $\kappa_{0}=6m$.
\medskip

\item $H_{1}^{'}=\{x_{3}^{2},x_{2}^{3},x_{4}^{m+1},x_{3}x_{4}^{m},x_{2}x_{4}^{m}\}$ and $\mathfrak{B}_{1}= \{x_{2}x_{4}^{m}, x_{2}^{2}x_{4}^{m},x_{3}x_{4}^{m},x_{2}x_{3}x_{4}^{m}, x_{2}^{2}x_{3}x_{4}^{m}\}$.
 Hence $\kappa_{1}=6m+1$.
\medskip
 
\item $H_{2}^{'}=\{x_{3}^{2},x_{2}^{3},x_{4}^{m+1},x_{3}x_{4}^{m},x_{2}^{2}x_{4}^{m}\}$ and 
$\mathfrak{B}_{2}=\{ x_{2}^{2}x_{4}^{m},x_{3}x_{4}^{m},x_{2}x_{3}x_{4}^{m}, x_{2}^{2}x_{3}x_{4}^{m}\}$. 
Therefore $\kappa_{2}=6m+2$.
\medskip
 
\item $H_{3}^{'}=\{x_{3}^{2},x_{2}^{3},x_{4}^{m+1},x_{3}x_{4}^{m}\}$ and $\mathfrak{B}_{3}= \{x_{3}x_{4}^{m},x_{2}x_{3}x_{4}^{m}, x_{2}^{2}x_{3}x_{4}^{m}\}$. Hence $\kappa_{3}=6m+3$.
\medskip
 
\item $H_{4}^{'}=\{x_{3}^{2},x_{2}^{3},x_{4}^{m+1},x_{2}x_{3}x_{4}^{m}\}$ and $\mathfrak{B}_{4}=\{x_{2}x_{3}x_{4}^{m}, x_{2}^{2}x_{3}x_{4}^{m}\}$. Hence $\kappa_{4}=6m+4$.
\medskip
 
\item $H_{5}^{'}=\{x_{3}^{2},x_{2}^{3},x_{4}^{m+1},x_{2}^{2}x_{3}x_{4}^{m}\}$ and $\mathfrak{B}_{5}= \{ x_{2}^{2}x_{3}x_{4}^{m}\}$. Hence $\kappa_{5}=6m+5$. 
\end{enumerate} 

\noindent We now apply Lemma \ref{equal} to each case to prove that these indeed give us the 
minimal generating sets for the ideal $\mathcal{P}_{4}$ in various cases.  \qed
\bigskip
  
\section{Syzygies of $k[\Gamma_{4}]$}
\begin{lemma}\label{reg}
Suppose $a=6m$, and $\gcd(a,d)=1$. Then, the set 
$\{x_{3}^{2}-x_{1}^{2}x_{4},x_{2}^{3}-x_{1}^{3}x_{3},x_{1}^{4m+d}-x_{4}^{m}\}$ forms a regular 
sequence in $A$. 
\end{lemma}

\proof With respect to the lexicographic monomial order induced by 
$x_{2}>x_{3}>x_{1}>x_{4}$ on $A$, the leading terms of these polynomials 
are mutually coprime. Hence the set 
$\{x_{3}^{2}-x_{1}^{2}x_{4},x_{2}^{3}-x_{1}^{3}x_{3},x_{1}^{4m+d}-x_{4}^{m}\}$ 
forms a regular sequence.\qed
\medskip

\begin{corollary}\label{ci}
Suppose $a=6m$ and $\gcd(a,d)=1$, then the Koszul complex resolves 
$A/\mathfrak{p}_{4}$ and the Betti numbers are $\beta_{0}=1,\beta_{1}=3,\beta_{2}=3,\beta_{3}=1$. 
Hence the ring $A/\mathfrak{P}_{4}$ is complete intersection.
\end{corollary} 

\proof Proof follows from lemma \ref{reg}.\qed
\medskip

\begin{lemma}\label{regsequence}
Let $m,d$ be two positive integers; consider the polynomials $g_{1}=-x_{1}^{(4m+d+2)}+x_{2}x_{4}^{m}$, 
$g_{2}=x_{1}^{5}x_{4}-x_{2}^{3}x_{3}$ and $g_{3}=x_{1}^{(4m+d-1)}x_{2}^{2}-x_{3}x_{4}^{m}$. The 
set $\{g_{1},g_{2},g_{3}\}$ forms a regular sequence in $A=k[x_{1},x_{2},x_{3},x_{4}]$.
\end{lemma}

\proof Let $x_{3}>x_{1}>x_{2}>x_{4}$ induce the lexicographic monomial order on $A$. 
Then $\mathrm{Lt}(g_{1})=-x_{1}^{(4m+d+2)}$, $\mathrm{Lt}(g_{2})=-x_{2}^{3}x_{3}$ and 
$\mathrm{Lt}(g_{3})=-x_{3}x_{4}^{m}$. Since $\gcd(\mathrm{Lt}(g_{1}),\mathrm{Lt}(g_{2}))=1$, 
the set  $\{g_{1},g_{2}\}$ forms Gr\"{o}bner basis of $\mathfrak{G}$, with respect to the 
chosen monomial order and forms a regular sequence. Let $\mathfrak{G}=\langle g_{1},g_{2}\rangle$ 
and  $g_{3}h\in \mathfrak{G}$; we have to show that $h\in \mathfrak{G}$. After division we may 
assume that $\mathrm{Lt}(g_{1})\nmid \mathrm{Lt}(h)$ and $\mathrm{Lt}(g_{2})\nmid \mathrm{Lt}(h)$. 
Since $g_{3}h\in \mathfrak{G}$ and $\mathrm{Lt}(g_{1})\nmid \mathrm{Lt}(h)$, 
$\mathrm{Lt}(g_{2})\nmid \mathrm{Lt}(h)$, we have $x_{2}^{3}\mid \mathrm{Lt}(h)$ and 
$x_{3}\nmid \mathrm{Lt}(h)$. We write $h=m_{0}+\cdots +m_{r}$, where each $m_{i}$'s 
are monomials and $m_{0}>\cdots>m_{r}$, with respect to the chosen monomial order. 
Since $x_{3}>x_{1}>x_{2}>x_{4}$ is the lexicographic monomial order on $A$,  
$x_{3}\nmid m_{0}$ implies $x_{3}\nmid m_{i}$, for $1\leq i\leq r$. Suppose 
$x_{1}^{l_{i}}\mid m_{i}$ but $x_{1}^{l_{i}+1}\nmid m_{i}$, then $i<j$ implies 
$l_{j}\leq l_{i}<4m+d+2$. Let $m_{0}=x_{2}^{3}m_{0}^{'}$, then
$$(x_{1}^{(4m+d-1)}x_{2}^{2}-x_{3}x_{4}^{m})(x_{2}^{3}m_{0}^{'}+m_{1}+\cdots+m_{r})\in \mathfrak{G} .$$ After dividing by $g_{2}$ we get 
$$(x_{1}^{(4m+d-1)}x_{2}^{5}-x_{1}^{5}x_{4}^{m+1})m_{0}^{'}+(x_{1}^{(4m+d-1)}x_{2}^{2}-x_{3}x_{4}^{m})(m_{1}+\cdots+m_{r})\in \mathfrak{G}.$$ Then leading term of above polynomial is $-x_{3}x_{4}^{m}m_{1}$ and we can divide by $g_{2}$. Continuing this way we get 
$$(x_{1}^{(4m+d-1)}x_{2}^{5}-x_{1}^{5}x_{4}^{m+1})(m_{0}^{'}+\cdots+m_{r}^{'})\in \mathfrak{G},$$ 
where $m_{i}=x_{2}^{3}m_{i}^{'}$, for $0\leq i\leq r$. Notice that $m_{0}^{'}>\cdots>m_{r}^{'}$. 
If $m=d=1$, then $-x_{1}^{5}x_{4}^{2}m_{0}^{'}$, otherwise $x_{1}^{(4m+d-1)}x_{2}^{5}m_{0}^{'}$ 
is the leading term of the above polynomial.
\medskip

\noindent\textbf{Case 1.} 
Suppose $m=d=1$, then 
$$(x_{1}^{4}x_{2}^{5}-x_{1}^{5}x_{4}^{2})(m_{0}^{'}+\cdots+m_{r}^{'})\in \mathfrak{G}.$$ 
We have $\mathrm{Lt}(g_{1})=-x_{1}^{7}\mid -x_{1}^{5}x_{4}^{2}m_{0}^{'}$ 
(since $x_{3}\nmid m_{0}=\mathrm{Lt}(h)$), hence $x_{1}^{2}\mid m_{0}^{'}$. 
Let $m_{0}=x_{2}^{3}m_{0}^{'}=x_{1}^{2}x_{2}^{3}m_{0}^{''}$. After dividing by 
$g_{1}$ we get $$(x_{1}^{6}x_{2}^{5}-x_{2}x_{4}^{3})m_{0}^{''}+(x_{1}^{4}x_{2}-x_{1}^{5}x_{4}^{2})(m_{1}^{'}+\cdots+m_{r}^{'})\in \mathfrak{G}.$$ We continuously divide the above polynomial 
by $g_{1}$ and we get, 
$$(x_{1}^{6}x_{2}^{5}-x_{2}x_{4}^{3})(m_{0}^{''}+\cdots+m_{s}^{''} )+(x_{1}^{4}x_{2}-x_{1}^{5}x_{4}^{2})(m_{s+1}^{'}+\cdots+m_{r}^{'}) \in \mathfrak{G},$$ 
where $0\leq s\leq r$, and for $0\leq i\leq s$ we have $m_{i}^{'}=x_{1}^{2}m_{i}^{''}$ 
and the leading term is $x_{1}^{6}x_{2}^{5}m_{0}^{''}$. Therefore 
$x_{1}^{2}\nmid m_{i}^{'}$ for $s+1\leq i\leq r$. Hence 
$\mathrm{Lt}(g_{1})=-x_{1}^{7}\mid x_{1}^{6}x_{2}^{5}m_{0}^{''}$, 
therefore $x_{1}\mid m_{0}^{''}$. Let $m_{0}=x_{2}^{3}x_{1}^{3}m_{0}^{'''}$, 
then we have 
$$(x_{2}^{6}x_{4}-x_{4}^{3}x_{2}x_{1})m_{0}^{'''}+(x_{1}^{6}x_{2}^{5}-x_{2}x_{4}^{3})(m_{1}^{''}+\cdots+m_{s}^{''} )+(x_{1}^{4}x_{2}-x_{1}^{5}x_{4}^{2})(m_{s+1}^{'}+\cdots+m_{r}^{'}) \in \mathfrak{G}.$$ 
Again we continuously divide by $g_{1}$ and for $0\leq s\leq r$ we get 
$$(x_{2}^{6}x_{4}-x_{4}^{3}x_{2}x_{1})(m_{0}^{'''}+\cdots+m_{s}^{'''})+(x_{1}^{4}x_{2}-x_{1}^{5}x_{4}^{2})(m_{s+1}^{'}+\cdots+m_{r}^{'}) \in \mathfrak{G}.$$ 
If the leading term of above polynomial is $-x_{1}^{5}x_{4}^{2}m_{s+1}^{'}$, then 
$\mathrm{Lt}(g_{1})\nmid -x_{1}^{5}x_{4}^{2}m_{s+1}^{'} $ (as $x_{1}^{2}\nmid m_{s+1}^{'}$ ). 
Therefore the leading term of above polynomial is $-x_{4}^{3}x_{2}x_{1}m_{0}^{'''} $. 
Thus, we have $\mathrm{Lt}(g_{1})=-x_{1}^{7}\mid -x_{4}^{3}x_{2}x_{1}m_{0}^{'''}$, 
hence $x_{1}^{6}\mid m_{0}^{'''}$, which implies that $\mathrm{Lt}(g_{1})\mid \mathrm{Lt}(h)$ - 
a contradiction.
\medskip

\noindent\textbf{Case 2.} If $m$ or $d$ is greater than $1$, then 
$x_{1}^{(4m+d-1)}x_{2}^{5}m_{0}^{'}$ is the leading term of 
$(x_{1}^{(4m+d-1)}x_{2}^{5}-x_{1}^{5}x_{4}^{m+1})(m_{0}^{'}+\cdots+m_{r}^{'})$. 
Therefore $\mathrm{Lt}(g_{1})=-x_{1}^{4m+d+2}\mid  x_{1}^{(4m+d-1)}x_{2}^{5}m_{0}^{'}$, 
hence $x_{1}^{3}\mid m_{0}^{'}$. After dividing by $g_{1}$ we get 
$$(x_{2}^{6}x_{4}^{m}-x_{1}^{8}x_{4}^{m+1})m_{0}^{''}+
(x_{1}^{(4m+d-1)}x_{2}^{5}-x_{1}^{5}x_{4}^{m+1})(m_{1}^{'}+\cdots+m_{r}^{'})\in \mathfrak{G}.$$
We proceed along the same line of argument as in Case 1. The variable  
$x_{3}$ is not present in the polynomial in each step, the leading term 
is always divisible by $\mathrm{Lt}(g_{1})$ and after finite steps we get 
$\mathrm{Lt}(g_{1})\mid \mathrm{Lt}(h)$ - a contradiction.\qed
\medskip

\begin{proposition} Suppose $a=6$ and $d=1$. Then the complex, 
$$0\longrightarrow A^{2}\stackrel{\phi_{3}}\longrightarrow A^{6}\stackrel{\phi_{2}}\longrightarrow A^{5}\stackrel{\phi_{1}}\longrightarrow A\longrightarrow A/\mathfrak{p}_{4}\longrightarrow 0 $$ 
is a minimal graded free resolution of $A/\mathfrak{p}_{4}$, where the maps $\phi_{i}$ are given by
$$\phi_{1}=(f_{1},f_{2},f_{3},f_{4},f_{5}),$$ with 
$f_{1}=-x_{2}^{3}+x_{1}^{3}x_{3}$, $f_{2}=-x_{3}^{2}+x_{1}^{2}x_{4}$, 
$f_{3}=x_{1}^{7}-x_{2}x_{4}$, $f_{4}=x_{1}^{4}x_{2}^{2}-x_{3}x_{4}$, 
$f_{5}=x_{1}^{2}x_{2}^{2}x_{3}-x_{4}^{2}$; 

$$\phi_{2}=\begin{pmatrix}
x_{1}^{4}& x_{4} &0& x_{1}^{2}x_{3}& 0 & x_{3}^{2}\\
0 & 0&x_{4}& x_{1}^{5}& x_{1}^{2}x_{2}^{2}&x_{1}^{3}x_{3}-x_{2}^{3}\\
-x_{3}&-x_{2}^{2}&0&-x_{4}&0&-x_{1}^{2}x_{2}^{2}\\
x_{2}&x_{1}^{3}&-x_{3}&0&-x_{4}&x_{1}^{5}\\
0&0&x_{1}^{2}&x_{2}&x_{3}&0
\end{pmatrix} $$ and 
$$\phi_{3}=\begin{pmatrix}
 x_{4}&x_{2}^{2}x_{3}\\
 -x_{1}^{4} &-x_{3}^{2}\\
 0& -x_{1}^{3}x_{3}+x_{2}^{3}\\
 -x_{3}& -x_{1}^{2}x_{2}^{2}\\
 x_{2}&x_{1}^{5}\\
 x_{1}^{2}&x_{4}
\end{pmatrix}.$$
\end{proposition} 

\proof We use the Buchsbaum-Eisenbud acyclicity theorem (see in \cite{bh}). 
It is easy to show that $\mathrm{grade}(I_{4}(\phi_{2}), A)\geq 2$. We take the minors, 
$$[1234 \mid 1235]=(x_{3}x_{4}-x_{1}^{4}x_{2}^{2})(x_{1}^{2}x_{2}^{2}x_{3}-x_{4}^{2}), 
[2345 \mid 1236]=x_{1}^{2}(x_{1}^{3}x_{3}-x_{2}^{3})^{2},$$ which 
have distinct irreducible factors, hence they form a regular sequence. 
We now consider the following minors,
$$[56 \mid 12]=-x_{1}^{7}+x_{2}x_{4}, [15\mid 12]=x_{1}^{5}x_{4}-x_{2}^{3}x_{3}, 
[46 \mid 12]=x_{1}^{4}x_{2}^{2}-x_{3}x_{4}^{m},$$ 
which form a regular sequence by Lemma \ref{regsequence}. \qed
\medskip

\begin{proposition}\label{syz1} Suppose $a=6m+1$ and either $m$ or $d$ is greater than $1$. 
Then the complex
$$0\longrightarrow A^{2}\stackrel{\phi_{3}}\longrightarrow A^{6}\stackrel{\phi_{2}}\longrightarrow A^{5}\stackrel{\phi_{1}}\longrightarrow A\longrightarrow A/\mathfrak{p}_{4}\longrightarrow 0$$ 
is a minimal graded free resolution of $A/\mathfrak{p}_{4}$, where the maps $\phi_{i}$ are 
given by
$$\phi_{1}=(f_{1},f_{2},f_{3},f_{4},f_{5}),$$ 
with 
$f_{1}=-x_{2}^{3}+x_{1}^{3}x_{3}$, $f_{2}=-x_{3}^{2}+x_{1}^{2}x_{4}$, 
$f_{3}=x_{1}^{(4m+d+2)}-x_{2}x_{4}^{m}$, $f_{4}=x_{1}^{(4m+d-1)}x_{2}^{2}-x_{3}x_{4}^{m}$, 
$f_{5}=x_{1}^{(4m+d-6)}x_{2}^{5}-x_{4}^{m+1}$;

$$\scriptsize{\phi_{2}=\begin{pmatrix}
x_{1}^{2}x_{4}-x_{3}^{2}& x_{1}^{(4m+d-1)} & x_{4}^{m}& x_{1}^{(4m+d-4)}x_{2}^{2} & x_{1}^{(4m+d-3)}x_{3}+x_{1}^{(4m+d-6)}x_{2}^{3} &x_{1}^{(4m+d-6)}x_{2}^{2}x_{3}\\
-x_{1}^{3}x_{3}+x_{2}^{3}& 0& 0& x_{4}^{m}& x_{1}^{4m+d}& x_{1}^{(4m+d-3)}x_{2}^{2}\\
0&-x_{3}&-x_{2}&0&-x_{4}&0\\
0&x_{2}&x_{1}^{3}&-x_{3}&0&-x_{4}\\
0&0&0&x_{1}^{2}&x_{2}&x_{3}
\end{pmatrix}}$$ 
and 

$$\scriptsize{\phi_{3}=\begin{pmatrix}
 x_{1}^{(4m+d-3)}&x_{4}^{m}\\
 -x_{4} &-x_{2}^{2}x_{3}\\
 0& -x_{1}^{2}x_{4}+x_{3}^{2}\\
 0& x_{1}^{3}x_{3}-x_{2}^{3}\\
 x_{3}&x_{1}^{2}x_{2}^{2}\\
 -x_{2}&-x_{1}^{5}
\end{pmatrix}.}$$
\end{proposition} 

\proof  We use the Buchsbaum-Eisenbud acyclicity theorem. It is easy to show that 
$\mathrm{grade}(I_{4}(\phi_{2}), A)\geq 2$. We take the minors 
$[1345 \mid 1246] = -x_{3}(x_{1}^{2}x_{4}-x_{3}^{2})^{2}$ and 
$[2345 \mid 2345] = (x_{4}^{m}x_{2}-x_{1}^{(4m+d+2)})(-x_{1}^{3}x_{3}+x_{2}^{3})$. These 
have distinct irreducible factors. Hence they form a regular sequence. 
Next we have to show that $\mathrm{grade}(I_{2}(\phi_{3}), A)\geq 3$. By Lemma \ref{regsequence}, 
the minors $[16\mid 12]=-x_{1}^{(4m+d+2)}+x_{2}x_{4}^{m}$, 
$[26\mid 12]=x_{1}^{5}x_{4}-x_{2}^{3}x_{3}$ and $[15\mid 12]=x_{1}^{(4m+d-1)}x_{2}^{2}-x_{3}x_{4}^{m}$ 
form a regular sequence. \qed
\medskip

\begin{proposition} Suppose $a=6m+2$. Then the complex, 
$$0\longrightarrow A^{2}\stackrel{\phi_{3}}\longrightarrow A^{6}\stackrel{\phi_{2}}\longrightarrow A^{5}\stackrel{\phi_{1}}\longrightarrow A\longrightarrow A/\mathfrak{p}_{4}\longrightarrow 0 $$ 
is a minimal graded free resolution of $A/\mathfrak{p}_{4}$, where the maps $\phi_{i}$ are given by 
$$\phi_{1}=(f_{1},f_{2},f_{3},f_{4},f_{5}),$$ with $f_{1}=-x_{2}^{3}+x_{1}^{3}x_{3}$, 
$f_{2}=-x_{3}^{2}+x_{1}^{2}x_{4}$, 
$f_{3}=x_{1}^{(4m+d+1)}x_{2}-x_{3}x_{4}^{m}$, 
$f_{4}=x_{1}^{(4m+d+4)}-x_{2}^{2}x_{4}^{m}$, 
$f_{5}=x_{1}^{(4m+d-4)}x_{2}^{4}-x_{4}^{m+1}$; 

$$\scriptsize{\phi_{2}=\begin{pmatrix}
x_{1}^{2}x_{4}-x_{3}^{2}& x_{4}^{m} & x_{1}^{(4m+d-2)}x_{2}& x_{1}^{(4m+d+1)} & x_{1}^{(4m+d-4)}x_{2}x_{3} &x_{1}^{(4m+d-1)}x_{3}-x_{1}^{(4m+d-4)}x_{2}^{3}\\
-x_{1}^{3}x_{3}+x_{2}^{3}& 0&  x_{4}^{m}&0& x_{1}^{(4m+d-1)}x_{2}& x_{1}^{(4m+d+2)}\\
0&x_{1}^{3}&-x_{3}&x_{2}^{2}&-x_{4}&0\\
0&-x_{2}&0&-x_{3}&0&-x_{4}\\
0&0&x_{1}^{2}&0&x_{3}&x_{2}^{2}
\end{pmatrix}}$$ 
and 
$$\scriptsize{\phi_{3}=\begin{pmatrix}
 x_{1}^{(4m+d-1)}&x_{4}^{m}\\
0 &-x_{1}^{2}x_{4}+x_{3}^{2}\\
 0& x_{1}^{3}x_{3}-x_{2}^{3}\\
 -x_{4}& -x_{2}x_{3}\\
 -x_{2}^{2}&-x_{1}^{5}\\
 x_{3}&x_{1}^{2}x_{2}
\end{pmatrix}.}$$
\end{proposition}

\proof The proof is similar to that of Proposition \ref{syz1}. We note that the following minors  
$[1345 \mid 1235]=x_{2}(x_{1}^{2}x_{4}-x_{3}^{2})^{2}$, 
$[2345 \mid 2345]=(x_{4}^{m}x_{3}-x_{1}^{(4m+d+1)}x_{2})(-x_{1}^{3}x_{3}+x_{2}^{3})$ 
in $I_{4}(\phi_{2})$ form a regular sequence. We then show that the minors belonging to 
$I_{2}(\phi_{3})$, given by $[15 \mid 12]=-x_{1}^{(4m+d+4)}+x_{2}x_{4}^{m}$, 
$[45 \mid 12]=x_{1}^{5}x_{4}-x_{2}^{3}x_{3}$ and $[16\mid 12]=x_{1}^{(4m+d+1)}x_{2}^{2}-x_{3}x_{4}^{m}$,  
form a regular sequence.\qed
\medskip

\begin{proposition} Suppose $a=6m+3$. Then the complex, 
$$0\longrightarrow A^{2}\stackrel{\phi_{3}}\longrightarrow A^{5}\stackrel{\phi_{2}}\longrightarrow A^{4}\stackrel{\phi_{1}}\longrightarrow A\longrightarrow A/\mathfrak{p}_{4}\longrightarrow 0 $$ is a 
minimal graded free resolution of $A/\mathfrak{p}_{4}$, where the maps $\phi_{i}$ are given by 
$$\phi_{1}=(f_{1},f_{2},f_{3},f_{4}),$$ with $f_{1}=-x_{2}^{3}+x_{1}^{3}x_{3}$, 
$f_{2}=-x_{3}^{2}+x_{1}^{2}x_{4}$, $f_{3}=x_{1}^{(4m+d+3)}-x_{3}x_{4}^{m}$, 
$f_{4}=x_{1}^{(4m+d-2)}x_{2}^{3}-x_{4}^{m+1}$; 

$$\phi_{2}=\begin{pmatrix}
x_{1}^{2}x_{4}-x_{3}^{2}&x_{1}^{(4m+d)} & x_{3}x_{4}^{m}& x_{1}^{(4m+d-2)}x_{3} & x_{1}^{(4m+d-2)}x_{2}^{3}-x_{4}^{m+1} \\
-x_{1}^{3}x_{3}+x_{2}^{3}& x_{4}^{m}&x_{1}^{3}x_{4}^{m}& x_{1}^{(4m+d+1)}&0\\
0&-x_{3}&-x_{2}^{3}&-x_{4}&0\\
0&x_{1}^{2}&x_{1}^{5}&x_{3}&-x_{1}^{3}x_{3}+x_{2}^{3}\\
\end{pmatrix} $$ 
and 
$$\scriptsize{\phi_{3}=\begin{pmatrix}
x_{4}^{m}& x_{1}^{(4m+d+1)}\\
-x_{2}^{3} &-x_{1}^{3}x_{4}\\
 x_{3}& x_{4}\\
 0& x_{1}^{3}x_{3}-x_{2}^{3}\\
 x_{1}^{2}&x_{3}\\
\end{pmatrix}.}$$
\end{proposition}

\proof The proof is similar to that of Proposition \ref{syz1}. We observe that the minors 
$[134\mid 124]=(x_{1}^{2}x_{4}-x_{3}^{2})^{2}$ and 
$[234\mid 123]=x_{1}^{2}(-x_{1}^{3}x_{3}+x_{2}^{3})^{2}$ belonging to 
$I_{3}(\phi_{2})$ form a regular sequence. We further observe that the minors 
$[15\mid 12]=-x_{1}^{(4m+d+3)}+x_{3}x_{4}^{m}$, 
$[13\mid 12]=x_{4}^{m+1}-x_{1}^{4m+d+1}x_{3}$ and 
$[34\mid 12]=x_{1}^{3}x_{3}^{2}-x_{2}^{3}x_{3}$, belonging to 
$I_{2}(\phi_{3})$, form a regular sequence.\qed
\medskip

\begin{proposition} Suppose $a=6m+4$. Then the complex, 
$$0\longrightarrow A^{3}\stackrel{\phi_{3}}\longrightarrow A^{6}\stackrel{\phi_{2}}\longrightarrow A^{4}\stackrel{\phi_{1}}\longrightarrow A\longrightarrow A/\mathfrak{p}_{4}\longrightarrow 0 $$ 
is minimal graded free resolution of $R/\mathcal{p}_{4}$, where the maps $\phi_{i}$ are given by
$$\phi_{1}=(f_{1},f_{2},f_{3},f_{4}),$$ with $f_{1}=-x_{2}^{3}+x_{1}^{3}x_{3}$, 
$f_{2}=-x_{3}^{2}+x_{1}^{2}x_{4}$, $f_{3}=x_{1}^{(4m+d+5)}-x_{2}x_{3}x_{4}^{m}$, 
$f_{4}=x_{1}^{(4m+d)}x_{2}^{2}-x_{4}^{m+1}$;

$$\scriptsize{\phi_{2}=\begin{pmatrix}
x_{1}^{2}x_{4}-x_{3}^{2}&  x_{1}^{(4m+d+2)} & x_{3}x_{4}^{m}& x_{1}^{(4m+d)}x_{3} & x_{1}^{(4m+d)}x_{2}^{2}-x_{4}^{m+1} &0\\
-x_{1}^{3}x_{3}+x_{2}^{3}& x_{4}^{m}x_{2}& x_{1}^{3}x_{4}^{m}&x_{1}^{(4m+d+3)}& 0&x_{1}^{(4m+d)}x_{2}^{2}-x_{4}^{m+1}\\
0&-x_{3}&-x_{2}^{2}&-x_{4}&0&0\\
0&x_{1}^{2}x_{2}&x_{1}^{5}&x_{2}x_{3}&-x_{1}^{3}x_{3}+x_{2}^{3}&-x_{1}^{2}x_{4}+x_{3}^{2}
\end{pmatrix} }$$ 
and 
$$\scriptsize{\phi_{3}=\begin{pmatrix}
x_{4}^{m}&0& x_{1}^{(4m+d)}\\
-x_{2}^{2} &0&-x_{4}\\
 x_{3}& x_{4}&0\\
 0& -x_{2}^{2}&x_{3}\\
 x_{1}^{2}&x_{3}&0\\
 0&x_{1}^{3}&-x_{2}
\end{pmatrix}.}$$
\end{proposition}

\proof The proof is similar to that of Proposition \ref{syz1}. We observe that the minors 
$[134\mid 124]=x_{2}(x_{1}^{2}x_{4}-x_{3}^{2})^{2}$ and 
$[234\mid 123]=x_{1}^{2}(-x_{1}^{3}x_{3}+x_{2}^{3})^{2}$ in $I_{3}(\phi_{2})$ form 
a regular sequence. We also observe that the minors 
$[124\mid 123]=x_{1}^{(4m+d)}x_{2}^{4}-x_{2}^{2}x_{4}^{m+1}$, 
$[246\mid 123]=x_{1}^{3}x_{2}^{2}x_{3}-x_{2}^{5}$ and 
$[256\mid 123]=x_{2}^{3}x_{3}-x_{1}^{5}x_{4}$, belonging to $I_{3}(\phi_{3})$,  
form a regular sequence. \qed
\medskip

\begin{proposition} Suppose $a=6m+5$. Then the complex, 
$$0\longrightarrow A^{3}\stackrel{\phi_{3}}\longrightarrow A^{6}\stackrel{\phi_{2}}\longrightarrow A^{4}\stackrel{\phi_{1}}\longrightarrow A\longrightarrow A/\mathcal{P}_{4}\longrightarrow 0 $$ is 
a minimal graded free resolution of $A/\mathcal{p}_{4}$, where the maps $\phi_{i}$ are given by
$$\phi_{1}=(f_{1},f_{2},f_{3},f_{4}),$$ with $f_{1}=-x_{2}^{3}+x_{1}^{3}x_{3}$, 
$f_{2}=-x_{3}^{2}+x_{1}^{2}x_{4}$, $f_{3}=x_{1}^{(4m+d+2)}x_{2}-x_{4}^{m+1}$, 
$f_{4}=x_{1}^{(4m+d+7)}-x_{2}^{2}x_{3}x_{4}^{m}$;

$$\scriptsize{\phi_{2}=\begin{pmatrix}
x_{1}^{2}x_{4}-x_{3}^{2}& x_{3}x_{4}^{m}& x_{1}^{(4m+d+4)}  & x_{1}^{(4m+d+2)}x_{2}-x_{4}^{m+1}&0 & x_{1}^{(4m+d+2)}x_{3}\\
-x_{1}^{3}x_{3}+x_{2}^{3}& x_{1}^{3}x_{4}^{m}&x_{4}^{m}x_{2}^{2}&0 &x_{1}^{(4m+d+2)}x_{2}-x_{4}^{m+1}& x_{1}^{(4m+d+5)}\\
0&x_{1}^{5}&x_{1}^{2}x_{2}^{2}&-x_{1}^{3}x_{3}+x_{2}^{3}&-x_{1}^{2}x_{4}+x_{3}^{2}&x_{2}^{2}x_{3}\\
0&-x_{2}&-x_{3}&0&0&-x_{4}
\end{pmatrix} }$$ 
and 
$$\scriptsize{\phi_{3}=\begin{pmatrix}
x_{4}^{m}&0& x_{1}^{(4m+d+2)}\\
x_{3}& x_{4}&0\\
-x_{2} &0&-x_{4}\\
 x_{1}^{2}&x_{3}&0\\
 0&x_{1}^{3}&-x_{2}^{2}\\
 0& -x_{2}&x_{3}\\
\end{pmatrix}.}$$
\end{proposition}

\proof The proof is similar to that of Proposition \ref{syz1}. We observe that the minors 
$[134\mid 125]=x_{2}(x_{1}^{2}x_{4}-x_{3}^{2})^{2}$ and 
$[234\mid 123]=x_{1}^{2}(-x_{1}^{3}x_{3}+x_{2}^{3})^{2} $
belonging to $I_{3}(\phi_{2})$ form a regular sequence. 
We further observe that the minors 
$[123\mid 123]=x_{1}^{(4m+d+2)}x_{2}x_{4}-x_{4}^{m+2}$, 
$[345\mid 123]=x_{2}^{3}x_{3}-x_{1}^{5}x_{4}$, 
$[456\mid 123]=x_{1}^{5}x_{3}-x_{1}^{2}x_{2}^{3}$,
belonging to $I_{3}(\phi_{3})$, form a regular sequence. \qed
\medskip

\begin{lemma}\label{stci}
The curve $k[\Gamma_{4}]$ is a set-theoretic complete intersection if 
$a\equiv i(\mathrm{mod}\, a)$, where $i\in\{0,3,4,5\}$.
\end{lemma}

\proof If $a=6m$, then $k[\Gamma_{4}]$ is a set-theoretic complete intersection by Corollary \ref{ci}. 
For the other cases it follows from theorem 5.3. in \cite{eto}.\qed
\bigskip

\section{Ap\'{e}ry table and Tangent Cone of $k[\Gamma_{4}]$}
Throughout this section we assume that the field $k$ is infinite.
\medskip

\begin{definition} Let $(R,m)$ be a Noetherian local ring and $I\subset R$ 
be an ideal of $R$. The fibre cone of $I$ is the ring 
$$F(I)=\displaystyle\bigoplus_{n\geq 0}\dfrac{I^{n}}{mI^{n}}\cong R[It]\otimes R/m.$$ 
Krull dimension of the ring $F(I)$ is called the \textit{analytic spread} of 
$I$, denoted by $\ell(I)$.
\end{definition}
\medskip

An ideal $J\subset I$ is called a \textit{reduction} of $I$ if there exists an 
integer $n>0$ such that $I^{n+1}=JI^{n}$. A reduction $J$ of $I$ is a \textit{minimal reduction} 
if $J$ is minimal with respect to inclusion among reductions of $I$. 
A minimal reduction always exists by \cite{nr}. It is well known that $J$ is a 
minimal reduction of $I$ if and only if $J$ is minimally generated by $\ell(I)$ 
number of elements, i.e, $\mu(J)=\ell(I)$. If  $J$ is a minimal reduction of $I$, 
then the least integer $r$ such that $I^{r+1}=JI^{r}$, is the reduction number 
of $I$ with respect to $J$, denoted by $r_{J}(I)$.
\medskip

We are particilarly interested in the semigroup ring $k[[\Gamma_{4}]]$, 
which is the coordinate ring of the algebroid monomial curve defined by 
the numerical semigroup $\Gamma_{4}$. Let  $a\geq 7$ and $d > 0$ be two integers, 
such that $\gcd(a,d)=1$. Let $R=k[[t^{a},t^{2a+d},t^{3a+3d},t^{4a+6d}]]$ and 
$\mathfrak{m}$ is the maximal ideal $\langle t^{a},t^{2a+d},t^{3a+3d},t^{4a+6d}\rangle$. 
Consider the principal ideal $I=\langle t^{a}\rangle \subset R$. The fibre cone of 
$I$ is the ring 
$$F(I)=\displaystyle\bigoplus_{n\geq 0}\dfrac{I^{n}}{\mathfrak{m}I^{n}}\cong R 
[It]\otimes R/\mathfrak{m}.$$ 
We note that here $\ell(I)=1$ and the tangent cone $G_{\mathfrak{m}}=\displaystyle\bigoplus_{n\geq 0}\dfrac{\mathfrak{m}^{n}}{\mathfrak{m}^{n+1}}$ is an $F(I)$-algebra. Moreover $F(I)\hookrightarrow G_{\mathfrak{m}} $ is a Noether normalisation (see \cite{cz1} and \cite{cz2}).
\medskip

Suppose $\Gamma$ be a numerical semigroup minimally generated by $a_{1}<\cdots 
<a_{e}$. Let $M=\Gamma\setminus\{0\}$ and for a positive integer $n$, we write 
$nM:=M+\cdots+M$ ($n$-copies). Let $\mathfrak{m}$ be the maximal ideal of the 
ring $k[[t^{a_{1}},\ldots t^{a_{e}}]]$. Then $(n+1)M=a+nM$ for all $n\geq r$ 
if and only if $r=r_{(t^{a_{1}})}(\mathfrak{m})$.
\medskip

Let $\mathrm{Ap}(\Gamma,a_{1})=\{0,\omega_{1},\ldots,\omega_{a_{1}-1}\}$. Now for 
each $n\geq 1$, let us define $\mathrm{Ap}(nM)=\{\omega_{n,0},\ldots\omega_{n,a_{1}-1}\}$ 
inductively. We define $\omega_{1,0}=a_{1}$ and $\omega_{1,i}=\omega_{i}$, for $1\leq i\leq a_{1}-1$. 
Then $\mathrm{Ap}(M)=\{a_{1},\omega_{1},\ldots,\omega_{a_{1}-1}\}$. Now we define 
$\omega_{n+1,i}=\omega_{n,i}$, if $\omega_{n,i}\in (n+1)M $, and $\omega_{n+1,i}=\omega_{n,i}+a_{1}$, 
otherwise. We note that $\omega_{n+1,i}=\omega_{n,i}+a_{1} $ for all $0\leq i\leq a_{1}-1$ and $n\leq r_{(t^{a_{1}})}(\mathfrak{m})$. Then, the Ap\'{e}ry table $\mathrm{AT}(\Gamma,a_{1})$ of $\Gamma$ is a table 
of size $(r_{(t^{a_{1}})}(\mathfrak{m})+1)\times a_{1}$, whose $(0,t)$ entry is $\omega_{t}$, 
for $0\leq t\leq {a_{1}-1}$ (we take $\omega_{0}=0$), and the $(s,t)$ entry is $\omega_{st}$, 
for $1\leq s\leq r_{(t^{a_{1}})}(\mathfrak{m})$ and $ 0\leq t\leq {a_{1}-1}$.
\medskip

Next we want to describe Apery table of $\Gamma_{4}$ and we need following Lemmas.
\medskip

\begin{lemma}\label{aperytable}
Elements of the Ap\'{e}ry set $\mathrm{Ap}(\Gamma_{4},a)$ have unique expressions.
\end{lemma}

\proof  We have $(4\mu_{i}+3\nu_{i}+2\xi_{i})a+id=\mu_{i}(4a+6d)+\nu_{i}(3a+3d)+\xi_{i}(2a+d)$, for $i\leq i\leq a-1$. Suppose for some $i\leq i\leq a-1$, $$(4\mu_{i}+3\nu_{i}+2\xi_{i})a+id=c_{1}(2a+d)+c_{2}(3a+3d)+c_{3}(4a+6d) .$$ Then
\begin{equation}\label{aptableeq1}
[(4c_{3}+3c_{2}+2c_{1})- (4\mu_{i}+3\nu_{i}+2\xi_{i})]a=[(6\mu_{i}+3\nu_{i}+\xi_{i})-(6c_{3}+3c_{2}+c_{1})]d.
\end{equation}
We have already shown in Theorem \ref{apery}, that, $4c_{3}+3c_{2}+2c_{1}\geq 4\mu_{i}+3\nu_{i}+2\xi_{i} $ and $ 6c_{3}+3c_{2}+c_{1}\geq 6\mu_{i}+3\nu_{i}+\xi_{i}$. From equation \ref{aptableeq1} we have, 
 \begin{equation}\label{aptableeq2}
 4c_{3}+3c_{2}+2c_{1}= 4\mu_{i}+3\nu_{i}+2\xi_{i}
 \end{equation}
 \begin{equation}\label{aptableeq3}
 6c_{3}+3c_{2}+c_{1}= 6\mu_{i}+3\nu_{i}+\xi_{i}.
 \end{equation}
 We eliminate $\mu_{i}$ and get,  
 \begin{equation}\label{aptableeq4}
 3(c_{2}-\nu_{i})= 4(\xi_{i}-c_{1}).
 \end{equation}
Let $c_{2}-\nu_{i}=4k$ and $\xi_{i}-c_{1}=3k$, for  $k\in \mathbb{Z}$. If 
$k>0$ then $\xi_{i}=3k+c_{1}$, which is impossible, since $0\leq \xi_{i}\leq 2$. 
If $k<0$ then $\nu_{i}=c_{2}-4k$, again a contradiction, since $0\leq \nu_{i}\leq 1$. 
Therefore $k=0$, hence $(\mu_{i},\nu_{i},\xi_{i})=(c_{3},c_{2},c_{1})$. \qed
\medskip

\begin{lemma}\label{red}
 Let $(\mu_{i},\nu_{i},\xi_{i})$ be the same as in Theorem \ref{apery}, 
 for $1\leq i\leq a-1$, and $(\mu_{0},\nu_{0},\xi_{0})=(0,0,0)$. Then 
 $\lfloor\frac{a}{6}\rfloor+2 =\max\{\mu_{i}+\nu_{i}+\xi_{i}\mid 1\leq i\leq a-1\}$, 
 where $\lfloor\centerdot\rfloor$ denotes the greatest integer function.
\end{lemma}
 
\proof Let $a=6\mu +q$, $0\leq q\leq 5$. We note that $i\leq 6\mu+4$ for all 
 $1\leq i\leq a-1$. Therefore, $\mu_{i}+\nu_{i}+\xi_{i}\leq \mu+2$, for all $1\leq i\leq a-1$. 
 On the other hand, $\mu_{a-q-1}+\nu_{a-q-1}+\xi_{a-q-1}=\mu+2=\lfloor\frac{a}{6}\rfloor+2$. \qed
 \medskip
 
\begin{corollary}  Let $\mathrm{AT}(\Gamma_{4})$ denote the Ap\'{e}ry table for 
 $\Gamma_{4}$. Then $\mathrm{AT}(\Gamma_{4})$ will be of order 
 $(\lfloor\frac{a}{6}\rfloor+3) \times a$. Let $\omega_{st}$ be the 
 $(s,t)$ entry of the table $\mathrm{AT}(\Gamma_{4})$. Then,  
 $\omega_{st} = (4\mu_{t}+3\nu_{t}+2\xi_{t})a+td$, if $0\leq s\leq \mu_{t}+\nu_{t}+\xi_{t}$ and 
 $0\leq t\leq a-1$. On the other hand, 
 $\omega_{st} = (3\mu_{t}+2\nu_{t}+\xi_{t}+s)a+td$, if $\mu_{t}+\nu_{t}+\xi_{t}<s\leq\lfloor\frac{a}{6}\rfloor+2$ 
 and $0\leq t\leq a-1$. Hence the reduction number of $r_{\mathfrak{I}}(\mathfrak{m})$ is $\lfloor\frac{a}{6}\rfloor+2$.
\end{corollary} 
 
\proof Proof follows from Lemmas \ref{aperytable} and \ref{red}.\qed 
\medskip
  
\noindent\textbf{Remark.} Minimal generating set of the defining ideal can be found abstractly in \cite{rs1}, when elements of Apery set has unique representation. But here we have written explicitly.
\medskip

\begin{lemma}\label{power}
Let $(\mu_{i},\nu_{i},\xi_{i})$ be the same as in Lemma \ref{red}, for $0\leq i\leq a-1$. 
Let $a=6\mu+q$, $\mu\geq 1$, $0\leq q\leq 5$. Let $t_{k}$ be the number of solutions 
of the equation $\mu_{i}+\nu_{i}+\xi_{i}=k$, for $0\leq k\leq \mu+2$. Then
\begin{eqnarray*}
t_{k} & = & \begin{cases}
1 & \mbox{if} \quad k=0,\\
3 & \mbox{if} \quad k=1,\\
\lfloor\frac{q}{2}\rfloor+2 & \mbox{if} \quad k=2 \, \mbox{and}\, \mu=1,\\
5 & \mbox{if} \quad k=2 \quad \mbox{and} \quad \mu\geq 2,\\
6 & \mbox{if} \quad 3\leq k\leq\mu,\\ 
\lfloor\frac{q}{2}\rfloor+3 & \mbox{if} \quad k=\mu+1 \, \mbox{and} \, \mu\geq 2,\\
1 & \mbox{if} \quad k=\mu+2 \, \mbox{and} \, q\in\{0,1,2\},\\
2 & \mbox{if} \quad k=\mu+2 \, \mbox{and} \, q\in\{3,4\},\\
3 & \mbox{if} \quad k=\mu+2 \, \mbox{and} \, q=5.
\end{cases}
\end{eqnarray*}
\end{lemma}

\proof For each of the following cases we write the set of solutions.
\begin{itemize}
\item[Case 1.] If $k=0$ then $(\mu_{i},\nu_{i},\xi_{i})=(0,0,0)$ is the only solution.

\item[Case 2.] If $k=1$ then $\{(1,0,0),(0,1,0),(0,0,1)\}$ is the set of solutions.

\item[Case 3.] If $k=2$ and $\mu\geq 2$, then $\{(2,0,0),(1,1,0),(1,0,1),(0,1,1),(0,0,2)\}$ is the set of solutions. 

\item[Case 4.] If $3\leq k\leq \mu$, 
then $\{(k,0,0),(k-1,1,0),(k-1,0,1),(k-2,1,1),(k-2,0,2),(k-3,1,2)\}$ is the set of solutions.

\item[Case 5.] If $k=\mu+1$ and $\mu\geq 2$, then, 
\begin{enumerate}
\item[(i)] if $q\in\{0,1\}$ then $\{(\mu-1,1,1),(\mu-1,0,2),(\mu-2,1,2)\}$ is the set of solutions; 
\item[(ii)] if $q\in\{2,3\}$ then $\{(\mu,0,1),(\mu-1,1,1),(\mu-1,0,2),(\mu-2,1,2)\}$ is the set of solutions;
\item[(iii)] if $q\in\{4,5\}$ then $\{(\mu,0,1),(\mu,1,0),(\mu-1,1,1),(\mu-1,0,2),(\mu-2,1,2)\}$ is the set of solutions.  
\end{enumerate} 

\item[Case 6.] If $k=\mu+2$, then, 
\begin{enumerate}
\item[(i)] if $q\in\{0,1,2\}$ then $\{(\mu-1,1,2)\}$ is the set of solutions; 
\item[(ii)] if $q\in\{3,4\}$ then $\{(\mu,0,2),(\mu-1,1,2)\}$ is the set of solutions;
\item[(iii)] if $q=5$ then $\{(\mu,1,1),(\mu,0,2),(\mu-1,1,2)\}$ is the set of solutions.  
\end{enumerate}
\end{itemize} 
For the case $\mu=1$ and $k=2$, it is easy to calculate 
(see example \ref{atexample}).
 \bigskip
 
We take some definitions from \cite{cz2}. Let $W =\{a_{0},\ldots,a_{n}\}$ be a set of integers. We call it a \textit{ladder} if $a_{0}\leq\ldots\leq a_{n}$. Given a ladder, we say that a subset $L=\{a_{i},\ldots,a_{i+k}\}$, with $k\geq 1$, is a \textit{landing} of length $k$ if $a_{i-1}<a_{i}=\cdots=a_{i+k}<a_{i+k+1}$ (where $a_{-1}= -\infty$ and $a_{n+1}=\infty$). In this case, $s(L)=i$ and $e(L)=i+k$. A landing $L$ is said to be a \textit{true landing} if $s(L)\geq 1$. Given two landings $L$ and $L^{'}$, we set $L<L^{'}$ if $s(L)<s(L^{'})$. Let $p(W)+1$ be the number of landings and assume that $L_{0}<\cdots<L_{p(W)}$ are the distinct landings. Then we define the following numbers:
$s_{j}(W)=s(L_{j})$, $e_{j}(W)=e(L_{j})$, for each $0\leq j\leq p(W)$;
$c_{j}(W)=s_{j}(W)-e_{j-1}(W)$, for each $0\leq j\leq p(W)$.
\medskip

Suppose $\Gamma$ be a numerical semigroup minimally 
generated by $a_{1}<\cdots <a_{e}$ and $\mathfrak{m}$ be the maximal ideal of $k[[t^{a_{1}},\ldots t^{a_{e}}]]$. Let $r= r_{(t^{a_{1}})}(\mathfrak{m})$,  $M=\Gamma\setminus\{0\}$ and 
$\mathrm{Ap}(nM)=\{\omega_{n,0},\ldots\omega_{n,a_{1}-1}\}$ for $0\leq n \leq r$. For every $1\leq i\leq a_{1}-1$, consider the ladder of the values $W^{i}=\{\omega_{n,i}\}_{0\leq n\leq r}$ and define the following integers:
\begin{enumerate}[(i)]
\item $p_{i}=p(W^{i})$
\item $d_{i}=e_{p_{i}}(W^{i})$
\item $b_{j}^{i}=e_{j-1}(W^{i})$ and 
$c_{j}^{i}=c_{j}(W^{i})$, for $1\leq j\leq p_{i}$.
\end{enumerate}
\medskip

\begin{theorem}\textbf{(Cortadellas, Zarzuela.)}\label{tangentcone} With the above notations, 
$$G_{m}\cong F\oplus\displaystyle\bigoplus_{i=1}^{a_{1}-1}\left(F(-d_{i})\displaystyle \bigoplus_{j=1}^{p_{i}}\dfrac{F}{(({t^{a_{1}})^{*})^{c_{j}^{i}}}F}(-b_{j}^{i})\right),$$
where $G_{m}$ is the tangent cone of $\Gamma$ and $F=F((t^{a_{1}}))$ is the fiber cone.
\end{theorem}

\proof See Theorem 2.3 in \cite{cz2}.\qed
\medskip

\begin{corollary}\label{exptangent} The tangent cone $G_{\mathfrak{m}}$ of $\Gamma_{4}$ is a free 
$F(\mathfrak{I})$-module. Moreover 
$$ G_{\mathfrak{m}}= \displaystyle\bigoplus_{k=0}^{\lfloor \frac{a}{6}\rfloor+2}(F(\mathfrak{I})(-k))^{t_{k}},$$ 
where $t_{k}$'s are given in Lemma \ref{power}. 
\end{corollary}  

\proof Proof follows from corollary \ref{aperytable} and \ref{tangentcone}.\qed
\medskip

\begin{corollary} The following properties hold good 
for the tangent cone $G_{\mathfrak{m}}$ of $\Gamma_{4}$ :
\begin{enumerate}
\item[(i)] $G_{\mathfrak{m}}$ is Cohen-Macaulay;
\item[(ii)] $G_{\mathfrak{m}}$ is not Gorenstein;
\item[(iii)] $G_{\mathfrak{m}}$ is Buchsbaum.
\end{enumerate}
\end{corollary}

\proof $(i)$ and $(ii)$ easily follow from the fact that  $G_{\mathfrak{m}}$ is a free $F(\mathfrak{I})$-module (see section 4 in \cite{cz2}). For proving $(ii)$, we use 
Theorem 20 in \cite{cz1}. Here we observe that if 
$G_{\mathfrak{m}}$ is  Gorenstein then $\mathfrak{m}^{n}\cap(\mathfrak{m}^{n+2}:\mathfrak{I})=\mathfrak{m}^{n+1} $, 
for $1\leq n\leq r_{\mathfrak{I}}(\mathfrak{m})$. Now 
$\mathfrak{m}^{n}\cap(\mathfrak{m}^{n+2}:\mathfrak{I})=\mathfrak{m}^{n+1} $, for all $1\leq n\leq r_{\mathfrak{I}}(\mathfrak{m})$ if and only if $nM_{4}\cap (n+2)M_{4}-a=(n+1)M_{4}$, for all $1\leq n\leq r_{\mathfrak{I}}(\mathfrak{m})$, where $M_{4}=\Gamma_{4}\setminus \{0\}$. Which is impossible, since $ (n+1)a\notin nM_{4}$. \qed
\medskip

\begin{corollary} Let $HG_{\mathfrak{m}}(x)$ be the Hilbert series of $G_{\mathfrak{m}}$. Then $$HG_{\mathfrak{m}}(x)=\displaystyle\left(\sum_{k=0}^{\lfloor \frac{a}{6}\rfloor+2} t_{k}x^{k}\right)/(1-x).$$
Where $t_{k}$'s are given in Lemma \ref{power}. 
\end{corollary}
\proof Follows from Corollay \ref{exptangent}.\qed
\medskip

\begin{example}\label{atexample}
Let us consider an example where $a=11$ and $d=24$. Hence $\Gamma_{4}=\langle 11,46,105,188\rangle$. Here $d\equiv 2(\mathrm{mod}\, a)$ and we have $\mathrm{Ap}(\Gamma_{4},a)=\{(4\mu_{i}+3\nu_{i}+2\xi_{i})a+id \mid 1\leq i\leq a-1 \}\cup \{0\}$, where $(\mu_{i},\nu_{i},\xi_{i})$ are same as in \ref{red}. Let $\omega_{i}=(4\mu_{i}+3\nu_{i}+2\xi_{i})a+id$ for $0\leq i\leq a-1$; the values are given in the table below:
\begin{center}
\begin{tabular}{|c|c|c|c|c|c|c|c|c|c|c|c|c|}
\hline
 $i$ & 0 & 1 & 2 & 3 & 4 & 5 & 6 & 7 & 8 & 9 & 10\\ 
\hline
$\xi_{i}$ & 0 & 1 & 2 & 0 & 1 & 2 & 0 & 1 & 2 & 0 & 1\\ 
\hline
$\nu_{i}$ & 0 & 0 & 0 & 1 & 1 & 1 & 0 & 0 & 0 & 1 & 1\\ 
\hline
$\mu_{i}$ & 0 & 0 & 0 & 0 & 0 & 0 & 1 & 1 & 1 & 1 & 1\\ 
\hline
$\omega_{i}$& 0 & 46 & 92 & 105 & 151 & 197 & 188 & 234 & 280 & 293 & 339\\
\hline
\end{tabular}
\end{center}
Let $M_{4}=\Gamma_{4}\setminus \{0\}$, then, we have,
\begin{center}
\begin{tabular}{|c|c|c|c|c|c|c|c|c|c|c|c|c|}
\hline
Ap($\Gamma_{4}$)& 0 & 46 & 92 & 105 & 151 & 197 & 188 & 234 & 280 & 293 & 339\\ 
\hline
Ap$(M_{4})$ & 11 & 46 & 92 & 105 & 151 & 197 & 188 & 234 & 280 & 293 & 339\\ 
\hline
Ap$(2M_{4})$ & 22 & 57 & 92 & 116 & 151 & 197 & 199 & 234 & 280 & 293 & 339\\ 
\hline
Ap$(3M_{4})$ & 33 & 68 & 103 & 127 & 162 & 197 & 210 & 245 & 280 & 304 & 339\\ 
\hline
\end{tabular}
\end{center}
\end{example}
\noindent From the Ap\'{e}ry table we get 
$G_{\mathfrak{m}}=F\oplus F(-1)^{3}\oplus F(-2)^{4}\oplus F(-3)^{3} $, where $F=F(t^{a})$, the fiber cone of $(t^{a})$. Therefore we have the Hilbert series
$$HG_{\mathfrak{m}}(x)=\dfrac{1+3x+4x^{2}+3x^{3}}{1-x}.$$
\bigskip

\section{Ap\'{e}ry set, Ap\'{e}ry table and the tangent cone of $k[\mathfrak{S}_{n+2}]$}
Let $a,d ,r,h$ be  positive integers with $\gcd(a,d)=\gcd(a,r)=1$ and $d>hn(r-1)$. 
Suppose $a_{0}=a$ and $a_{k+1}=ha+r^{k}d$, for $0\leq k\leq n$.
Let $\mathfrak{S}_{n+2}=\langle \{a_{0},a_{1},\ldots,a_{n+1}\}\rangle$ be the numerical semigroup with embedding dimension 
$n+2$, such that $\{a_{0}, a_{1},\ldots, a_{n+1}\}$ 
form a minimal system of generators for $\mathfrak{S}_{n+2}$.
\medskip

\begin{definition}
Let $m,r,n$ be positive integers and $m=\displaystyle\sum_{k=0}^{n}\alpha_{k}r^{k}$, where $0\leq \alpha_{i}\leq r-1$ for $i\in\{0,\ldots,n-1\}$. Then the expression $m=\displaystyle\sum_{k=0}^{n}\alpha_{k}r^{k}$ is called the $r$-adic representation of $m$ upto order $n$.
\end{definition}
\medskip

\begin{lemma}\label{r-adic}
Let $m$ and $r$ be two positive integers and $m=\displaystyle\sum_{k=0}^{n}\alpha_{k}r^{k}$ be the $r$-adic representation of $m$  upto order $n$. Then for any expression $m=\displaystyle\sum_{k=0}^{n}\beta_{k}r^{k}$, we have 
$$\displaystyle\sum_{k=0}^{n}\alpha_{k}\leq \displaystyle\sum_{k=0}^{n}\beta_{k}.$$ Moreover, 
$\displaystyle\sum_{k=0}^{n}\alpha_{k}< \displaystyle\sum_{k=0}^{n}\beta_{k}$, if $\displaystyle\sum_{k=0}^{n}\beta_{k}r^{k}$ 
is not an $r$-adic representation of $m$ upto order $n$. 
\end{lemma}

\proof  We proceed by induction on $n$. If $n=0$ then it follows trivially. At first we claim that $\beta_{n}\leq \alpha_{n}$. 
If not, then $\alpha_{n}+1\leq\beta_{n}$, hence  $(\alpha_{n}+1)r^{n}\leq\beta_{n}r^{n}$. Now 
$$m=\displaystyle\sum_{k=0}^{n}\alpha_{k}r^{k}\leq \displaystyle\sum_{k=0}^{n-1}(r-1)r^{k} +\alpha_{n}r^{n}= (r^{n}-1)+\alpha_{n}r^{n} <(\alpha_{n}+1)r^{n}\leq \beta_{n}r^{n},$$ which is a contradiction. Let $\alpha_{n}=t+\beta_{n}$, where $t\geq 0$. Again $\displaystyle\sum_{k=0}^{n-1}\alpha_{k}r^{k}+(t+\beta_{n})r^{n}=\displaystyle\sum_{k=0}^{n}\beta_{k}r^{k}$, therefore $\displaystyle\sum_{k=0}^{n-2}\alpha_{k}r^{k}+(tr+\alpha_{n-1})r^{n-1}=\displaystyle\sum_{k=0}^{n-1}\beta_{k}r^{k}$. By the induction hypothesis, $\displaystyle\sum_{k=0}^{n-2}\alpha_{k}+(tr+\alpha_{n-1})\leq \displaystyle\sum_{k=0}^{n-1}\beta_{k}$. Hence, we have,
\begin{align*}
\displaystyle\sum_{k=0}^{n}\alpha_{k}
&=\displaystyle\sum_{k=0}^{n-1}\alpha_{k}+t+\beta_{n}\\
&\leq \displaystyle\sum_{k=0}^{n-2}\alpha_{k}+(tr+\alpha_{n-1})+\beta_{n}\\
&\leq \displaystyle\sum_{k=0}^{n}\beta_{k}. \hspace*{3.5in} \qed
\end{align*} 
\medskip

\begin{theorem}\label{apery}
Let for each $i\in \{1,\ldots,a-1\}$, $i=\displaystyle\sum_{k=0}^{n}a_{ki}r^{k}$ be the $r$-adic representation of $i$ upto order $n$. Suppose $\ell_{i}=\displaystyle\sum_{k=0}^{n}a_{ki}$, for $1\leq i\leq a-1$. Then
 $\mathrm{Ap}(\mathfrak{S}_{n+2}, a)=\{\ell_{i}ha+id\mid 1\leq i\leq a-1\}\cup\{0\}$. 
 \end{theorem}
 
 \proof  Let $T=\{\ell_{i}ha+id\mid 1\leq i\leq a-1\}$. At first we note that $$\ell_{i}ha+id=\displaystyle\sum_{k=0}^{n}a_{ki}(ha+r^{k}d),\quad  1\leq i\leq a-1 .$$  Therefore $T\subset \mathfrak{S}_{n+2}$. Suppose $s\in \mathrm{Ap}(\mathfrak{S}_{n+2}, a)\setminus\{0\}$, with $s\equiv id(\mathrm{mod}a)$. Let $s=\displaystyle\sum_{k=0}^{n}c_{k+1}(ha+r^{k}d)$, then 
 $\gcd(a,d)=1$ forces that $\displaystyle\sum_{k=0}^{n}c_{k+1}r^{k}\equiv i(\mathrm{mod}a)$. Therefore $\displaystyle\sum_{k=0}^{n}c_{k+1}r^{k}= i+pa$, and we have $s=\displaystyle\left(\sum_{k=0}^{n}c_{k+1}\right)ha+(i+pa)d$. 
\medskip
 
If $p>0$ then 
 \begin{align*}
 s&=\displaystyle\left(\sum_{k=0}^{n}c_{k+1}\right)ha+(i+pa)d\\
 &\geq\displaystyle\left(\sum_{k=0}^{n}c_{k+1}+n(r-1)\right)ha+id \\
 &> nh(r-1)a+id\quad (\mathrm{as}\,\, s>0\,\,\mathrm{implies}\,\, \sum_{k=0}^{n}c_{k+1}>0)\\
 &\geq \ell_{i}+id.
 \end{align*}
This gives a contradiction as $s\in \mathrm{Ap}(\mathfrak{S}_{n+2}, a) $ and $s\equiv \ell_{i}+id(\mathrm{mod} a) $. If $p=0$,  then by Lemma \ref{r-adic}, we have $\ell_{i}\leq \displaystyle\sum_{k=0}^{n}c_{k+1} $. Therefore $s\geq \ell_{i}+id $. Now $s\in \mathrm{Ap}(\mathfrak{S}_{n+2}, a) $ and $s\equiv \ell_{i}+id(\mathrm{mod} a) $, therefore we have $s=\ell_{i}+id$, hence $s\in T$. \qed
\medskip

\begin{lemma}\label{unique}
Every element of $\mathrm{Ap}(\mathfrak{S}_{n+2})$ has a unique expression.
\end{lemma} 

\proof Let $$\omega(i)=\ell_{i}ha+id=\displaystyle\sum_{k=0}^{n}c_{k+1}(ha+r^{k}d)=(\displaystyle\sum_{k=0}^{n}c_{k+1})ha+(\displaystyle\sum_{k=0}^{n}c_{k+1}r^{k})d,$$ for $1\leq i\leq a-1$, where $\ell_{i}$'s are the same as in Theorem 
\ref{apery}.  Therefore, $\displaystyle\sum_{k=0}^{n}c_{k+1}r^{k}\equiv i(\mathrm{mod}a)$, hence $\displaystyle\sum_{k=0}^{n}c_{k+1}r^{k}= i+pa $ for some $p\geq 0$. If $p>0$ then,
 \begin{align*}
 \omega(i)&=\displaystyle\left(\sum_{k=0}^{n}c_{k+1}\right)ha+(i+pa)d\\
 &\geq\displaystyle\left(\sum_{k=0}^{n}c_{k+1}+n(r-1)\right)ha+id \\
 &> nh(r-1)a+id\quad (\mathrm{as}\,\, \omega(i)>0\,\,\mathrm{implies}\,\, \sum_{k=0}^{n}c_{k+1}>0)\\
 &\geq \ell_{i}+id.
 \end{align*}
This gives a contradiction as $\omega(i)\in \mathrm{Ap}(\mathfrak{S}_{n+2}, a) $. Therefore $\displaystyle\sum_{k=0}^{n}c_{k+1}r^{k}= i $.  If the expression $\displaystyle\sum_{k=0}^{n}c_{k+1}r^{k}$ is not an $r$-adic representation of $i$ upto order $n$, then $\displaystyle\sum_{k=0}^{n}c_{k+1}>\ell_{i} $ by lemma \ref{r-adic}, which is a contradiction. Therefore $\displaystyle\sum_{k=0}^{n}c_{k+1}r^{k}$ is 
an $r$-adic representation of $i$ upto order $n$ and by the uniqueness of $r$-adic representation, 
$\omega(i)$ upto order $n$ has unique expression for 
$1\leq i\leq a-1$.   \qed
\medskip

\begin{theorem}\label{aperytab}
Let $r=\max\{\ell_{i}\mid 1\leq i\leq a-1\} $, where $\ell_{0}=0$ and $\ell_{i}$'s are the same as in Theorem \ref{apery}, for $1\leq i\leq a-1$.  Let $\mathrm{AT}(\mathfrak{S}_{n+2},a)$ denote the Ap\'{e}ry table of $\mathfrak{S}_{n+2}$. Then $\mathrm{AT}(\mathfrak{S}_{n+2},a)$ will be of order $r \times a$. Let $\omega_{st}$ be the $(s,t)$ entry of the table $\mathrm{AT}(\mathfrak{S}_{n+2},a)$. Then 
 \begin{align*}
 \omega_{st} & = &  
 \begin{cases}
 \ell_{t}ha+td & \mbox{if} \quad 0\leq s\leq \ell_{t}, \, 0\leq t\leq a-1;\\
 \ell_{t}ha+td+(s-\ell_{t})a & \mbox{if} \quad 
 \ell_{t}<s\leq r, \, 0\leq t\leq a-1.
 \end{cases}
 \end{align*}
 \end{theorem}
 \proof Follows from Lemma \ref{unique}.\qed
\medskip

\begin{theorem} Let $k$ be an infinite field. 
The following properties hold for the tangent cone $G_{\mathfrak{m}}$ of $\mathfrak{S}_{n+2}$:
\begin{enumerate}
\item[(i)] $G_{\mathfrak{m}}$ is Cohen-Macaulay,
\item[(ii)] $G_{\mathfrak{m}}$ is not Gorenstein,
\item[(iii)] $G_{\mathfrak{m}}$ is Buchsbaum.
\end{enumerate}
\end{theorem}
\proof Follows from Theorem \ref{aperytab} and \cite{cz2}.\qed

\bibliographystyle{amsalpha}

\end{document}